\documentclass{amsart}

\usepackage{latexsym,amsmath,amstext,amsthm,pstricks,pstcol}
\usepackage[english]{babel}
\usepackage[T1]{fontenc}
\usepackage[dvips]{graphicx}
\usepackage{delarray}
\usepackage{stmaryrd}
\usepackage{vmargin}
\theoremstyle{definition}

\newtheorem{thm}{Theorem}[section]
\newtheorem{coro}[thm]{Corollary}
\newtheorem{lem}[thm]{Lemma}
\newtheorem{prop}[thm]{Proposition}
\newtheorem{definition}[thm]{Definition}
\newtheorem{rem}[thm]{Remark}

\title{Lipschitz extension of multiple Banach-valued functions in the sense of Almgren}

\author{Jordan Goblet}

\address{Jordan Goblet\\
Département de mathématique\\
Université catholique de Louvain\\
Chemin du cyclotron 2\\
1348 Louvain-La-Neuve (Belgium)}

\email{goblet@math.ucl.ac.be}
\date{\today}


\begin{document}

\maketitle

\begin{abstract} A
multiple-valued function $f:X\rightarrow {\bf Q}_Q(Y)$ is
essentially a rule assigning
 $Q$ unordered and non necessarily distinct elements of $Y$ to each element
  of $X$. We study the Lipschitz extension problem
  in this context by using two general Lipschitz extension theorems recently
 proved by U. Lang and T. Schlichenmaier. We prove that the pair
$\left(X,{\bf Q}_Q(Y)\right)$ has the Lipschitz extension property
if $Y$ is a Banach space and $X$ is a metric space with a finite
Nagata dimension. We also show that ${\bf Q}_{Q}(Y)$ is an
absolute Lipschitz retract if $Y$ is a finite algebraic
dimensional Banach space.
\end{abstract}

\section{Introduction}

A multiple-valued function in the sense of Almgren is a map of the
form $f:X\rightarrow {\bf Q}_{Q}(Y)$ where $X$ and $(Y,d)$ are
metric spaces. The particular target space is defined by
$${\bf Q}_{Q}(Y)=\left\{\sum_{i=1}^{Q}\llbracket
x_{i}\rrbracket\; :\; x_{i}\in Y \;\text{ for
}i=1,\hdots,Q\right\}$$ where $\llbracket x_{i}\rrbracket$ denotes
the Dirac measure at $x_i$ and it is equipped with the metric
$$\mathcal{S}\left(\sum_{i=1}^{Q}\llbracket
x_{i}\rrbracket,\sum_{j=1}^{Q}\llbracket
y_{j}\rrbracket\right)=\min
\left\{\max_{i=1,\hdots,Q}\;d\left(x_{i},y_{\sigma(i)}\right)\;
:\;\sigma\text{ is a permutation of }\{1,\hdots,Q\} \right\}.$$
Consequently a multiple-valued function $f:X\rightarrow {\bf
Q}_{Q}(Y)$ is essentially a rule assigning $Q$ unordered and non
necessarily distinct elements of $Y$ to each element of $X$. Such
maps are studied in complex analysis (see appendix 5 in \cite{H.
Whitney}). Indeed in complex function theory one often speaks of
the ``two-valued function $f(z)=z^{1/2}$ ''. This can be
considered as a function from $\mathbb{C}$ to  ${\bf
Q}_{2}(\mathbb{C})$.\par

 In his big regularity paper \cite{Frederick
J. Almgren}, F. J. Almgren introduced ${\bf
Q}_{Q}(\mathbb{R}^n)$-valued functions to tackle the problem of
estimating the size of the singular set of mass-minimizing
integral currents (see \cite{Frederick J. Almgren 2} for a
summary). Almgren's multiple-valued functions are a fundamental
tool for understanding geometric variational problems in
codimension higher than 1. The success of Almgren's regularity
theory raises the need of further studying multiple-valued
functions.\par

The Lipschitz extension problem asks for conditions on a pair of
metric spaces $X,Y$ such that every Lipschitz $Y$-valued function
defined on a subset of $X$ can be extended to all of $X$ with only a
bounded multiplicative loss in the Lipschitz constant. More
precisely the pair $(X,Y)$ is said to have the Lipschitz extension
property if there exists a constant $\lambda \geq 1$ such that for
every subset $A\subset X$, every Lipschitz function $f:A\rightarrow
Y$ can be extended to a Lipschitz function $F:X\rightarrow Y$ with
Lip$(F)\leq\lambda\;$Lip$(f)$. A metric space $Y$ is said to be an
absolute Lipschitz retract if for every metric space $X$, the pair
$(X,Y)$ has the Lipschitz extension property (see chapter 1 in
\cite{Y. Benjamini J. Lindenstrauss} for equivalent definitions).
This problem dates back to the work of Kirszbraun and Whitney in the
1930's, and has been extensively investigated in the last two
decades (see \cite{U. Lang T. Schlichenmaier} and \cite{lee naor}
for several recent breakthroughs).\par

In the present paper, we will be interested in the Lipschitz
extension of ${\bf Q}_Q(Y)$-valued functions when $Y$ is a Banach
space. In this context, an important remark is that a Lipschitz
${\bf Q}_Q(Y)$-valued function is much more than $Q$ glued Lipschitz
$Y$-valued functions. Indeed we noticed in \cite{J. Goblet} that the
following Lipschitz ${\bf Q}_{2}(\mathbb{R}^{2})$-valued function
$$f:  {\bf S}^{1}\subset\mathbb{R}^{2}  \rightarrow  {\bf Q}_{2}(\mathbb{R}^2): x=(\cos\theta,\sin\theta)
\mapsto f(x)=\left\llbracket
\left(\cos\frac{\theta}{2},\sin\frac{\theta}{2}\right)\right\rrbracket
+\left\llbracket
\left(-\cos\frac{\theta}{2},-\sin\frac{\theta}{2}\right)\right\rrbracket$$
doesn't split into two Lipschitz $\mathbb{R}^2$-valued branches.
Consequently the Lipschitz extension problems for ${\bf
Q}_{Q}(Y)$-valued functions and $Y$-valued functions are in general
two distinct problems.

In \cite{Frederick J. Almgren}, Almgren built an explicit
bilipschitz correspondence between ${\bf Q}_Q(\mathbb{R}^{n})$ and a
Lipschitz retract denoted $Q^{*}$ included in a Euclidean space.
This construction and McShane-Whitney's Theorem (see 2.10.44 in
\cite{Herbert Federer}) clearly imply that ${\bf
Q}_{Q}(\mathbb{R}^{n})$ is an absolute Lipschitz retract. For
multiple Banach-valued functions, such a bilipschitz correspondence
is not available.\par

In a recent paper \cite{U. Lang T. Schlichenmaier}, U. Lang and T.
Schlichenmaier obtained two general Lipschitz extension theorems
involving a Lipschitz connectedness assumption on the target space
and a bound on the
Nagata dimension denoted dim$_{N}$ below of either the source space or the target space:\\

\noindent {\bf Theorem 1.5 in \cite{U. Lang T. Schlichenmaier}}.
Suppose that $X,Y$ are metric spaces, dim$_{N}(X)\leq n<\infty$, and
$Y$ is complete. If $Y$ is Lipschitz $(n-1)$-connected, then
the pair $(X,Y)$ has the Lipschitz extension property.\\

\noindent {\bf Corollary 1.8 in \cite{U. Lang T. Schlichenmaier}}.
Suppose that $Y$ is a metric space with dim$_{N}(Y)\leq n<\infty$.
Then Y is an absolute Lipschitz retract if and only if $Y$ is
complete and Lipschitz $n$-connected.

\vspace{0.4cm} In Section \ref{section2} we prove that ${\bf
Q}_{Q}(Y)$ is complete in case $Y$ is. We recall what is meant by
Lipschitz connectedness and we prove that ${\bf Q}_{Q}(Y)$ enjoys
this property when $Y$ is a weakly convex geodesic space.\par
 In Section
\ref{nagata} we define the Nagata dimension and gather a number of
basic properties. We estimate the Nagata dimension of ${\bf
Q}_{Q}(Y)$ in accordance with the Nagata dimension of $Y$. We also
show that ${\bf Q}_{Q}(Y)$ has a finite Nagata dimension when $Y$ is
a finite algebraic dimensional Banach space.\par We finally combine
these results with Lang-Schlichenmaier's Theorems in order to prove
Theorem \ref{exto}
 and Theorem \ref{final}.\\

\begin{thm}\label{exto}
Suppose that $X$ is a metric space, dim$_{N}(X)<\infty$, and $Y$
is a complete weakly convex geodesic space. Then the pair
$\left(X,{\bf Q}_{Q}(Y)\right)$ has the Lipschitz extension
property. In particular, the pair $\left(X,{\bf Q}_{Q}(Y)\right)$
has the
Lipschitz extension property if $Y$ is a Banach space.\\
\end{thm}

\begin{thm}\label{final}
If $Y$ is a complete weakly convex geodesic space with a finite
Nagata dimension then ${\bf Q}_{Q}(Y)$ is an absolute Lipschitz
retract. In particular, ${\bf Q}_{Q}(Y)$ is an absolute Lipschitz
retract if $Y$ is a Banach space with a finite algebraic
dimension.
\end{thm}

\section{Completeness and Lipschitz connectedness of ${\bf
Q}_{Q}(Y)$}\label{section2}

For later use we note a simple fact related to the completeness
property.

\begin{lem}\label{complete} If $Y$
is a complete metric space, then ${\bf Q}_{Q}(Y)$ is complete.
\end{lem}
\vspace{0.4cm} \noindent {\bf Proof.}  Let
$(x^{n})_{n\in\mathbb{N}}$ be a Cauchy sequence in ${\bf
Q}_{Q}(Y)$. It is enough to show that we can extract a converging
subsequence. On the one hand, it is clear that we can extract a
subsequence $(x^{n_{l}})_{l\in\mathbb{N}}$ such that
$\mathcal{S}(x^{n_{l}},x^{n_{l+1}})<1/2^{l}$ for all
$l\in\mathbb{N}$. On the other hand, we can write
$x^{n_{l}}=\sum_{i=1}^{Q}\llbracket x_{i}^{n_{l}}\rrbracket$  with
$\mathcal{S}(x^{n_{l}},
x^{n_{l+1}})=\max_{i=1,\hdots,Q}\;d(x_{i}^{n_{l}},x_{i}^{n_{l+1}})$
for all $l\in\mathbb{N}$. We obtain that
$d(x_{i}^{n_{l}},x_{i}^{n_{l+1}})< 1/2^{l}$ hence
$(x_{i}^{n_{l}})_{l\in\mathbb{N}}$ is a Cauchy sequence in $Y$ for
$i=1,\hdots,Q$. Since $Y$ is complete, there exist
$x_{1},\hdots,x_{Q}\in Y$ such that $x_{i}^{n_{l}}\rightarrow
x_{i}$ as $l\rightarrow \infty$ for $i=1,\hdots,Q$. Consequently,
the subsequence $(x^{n_{l}})_{l\in\mathbb{N}}$ tends to
$\sum_{i=1}^{Q}\llbracket
x_{i}\rrbracket\in {\bf Q}_Q(Y)$ as $l\rightarrow \infty$.\qed\\

Recall that a topological space $Y$ is said to be $n$-connected,
for some integer $n\geq 0$, if for every $m\in\{0,1,\hdots,n\}$,
every continuous map from ${\bf S}^{m}$ into $Y$ admits a
continuous extension to ${\bf B}^{m+1}$. Accordingly, we call a
metric space $Y$ Lipschitz $n$-connected if there is a constant
$\lambda\geq 1$ such that for every $m\in\{0,1,\hdots,n\}$, every
Lipschitz map $f:{\bf S}^{m}\rightarrow Y$ possesses a Lipschitz
extension $F:{\bf B}^{m+1}\rightarrow Y$ with
Lip$(F)\leq\lambda\;$Lip$(f)$. Here ${\bf S}^{m}$ and ${\bf
B}^{m+1}$ denote the unit sphere and closed ball in
$\mathbb{R}^{m+1}$, equipped with the induced metric.\par

Let $(Y,d)$ be a metric space. A geodesic joining $x\in Y$ to
$y\in Y$ is a map $c_{xy}$ from a closed interval $[0,l]\subset
\mathbb{R}$ to $Y$ such that $c_{xy}(0)=x,\; c_{xy}(l)=y$ and
$d(c(t),c(t'))=|t-t'|$ for all $t,t'\in [0,l]$ (in particular,
$l=d(x,y)$). A geodesic bicombing on $Y$ is an assignment of a
geodesic $c_{xy}:[0,d(x,y)]\rightarrow Y$ from $x$ to $y$ for
every pair $(x,y)\in Y\times Y$. We call a geodesic bicombing
$\{c_{xy}\}$ on $Y$ $\gamma$-weakly convex, for some constant
$\gamma\geq 1$, if each pair of geodesics
$c_{xy}:[0,d(x,y)]\rightarrow Y$ and $c_{xz}:[0,d(x,z)]\rightarrow
 Y$ satisfy the inequality
\begin{equation}\label{waka}
d\left(c_{xy}(t\;d(x,y)),c_{xz}(t\;d(x,z))\right)\leq \gamma\;
t\;d(y,z)
\end{equation}
 for all $t\in [0,1]$. A metric space which admits a $\gamma$-weakly convex geodesic
bicombing is said to be a $\gamma$-weakly convex geodesic space.

\begin{rem}
On the one hand, the inequality (\ref{waka}) holds for $\gamma=1$
on every geodesic space with convex distance function (for the
unique geodesic bicombing), and on every normed vector space for
the linear geodesic bicombing. On the other hand, one readily
checks that a weakly convex geodesic space is Lipschitz
$n$-connected for all $n\in\mathbb{N}$.
\end{rem}

 Proposition \ref{connexe} generalizes a
construction described in Section 1.5 of \cite{Frederick J.
Almgren}. It proves that the metric space ${\bf Q}_{Q}(Y)$ is
Lipschitz $n$-connected for all $n\in\mathbb{N}$ if $Y$ is a weakly
convex geodesic space. The Lipschitz connectedness of ${\bf
Q}_{Q}(Y)$ is far from being obvious in this context since ${\bf
Q}_{Q}(Y)$ does not necessarily inherit from $Y$ a weakly convex
geodesic bicombing. Indeed ${\bf Q}_{2}(\mathbb{R}^2)$ does not
possess a weakly convex geodesic bicombing.

\begin{prop}\label{connexe}
Let $Y$ be a $\gamma$-weakly convex geodesic space. Every Lipschitz
multiple-valued function $f:{\bf S}^{m}\rightarrow {\bf Q}_Q(Y)$
extends to $F:{\bf B}^{m+1}\rightarrow {\bf Q}_Q(Y)$ with
 Lip$(F)\leq (\gamma+8Q-6)$ Lip$(f)$.
\end{prop}
\vspace{0.4cm}\noindent {\bf Proof.} Set
$D=2\;$Lip$(f)=\;$diam$({\bf S}^{m})\;$Lip($f$) and choose positive
integers $s,Q_{1},Q_{2},\hdots,Q_{s}$ and points $p(i,k)$ in $Y$ for
$k=1,\hdots,Q_{i},\;i=1,\hdots,s$ subject to the following
requirements:\\

\begin{enumerate}
\item
$f(1,0,0,\hdots,0)=\sum_{i=1}^{s}\sum_{k=1}^{Q_{i}}\llbracket
p(i,k)\rrbracket$;\\
\item \label{ii} If $i\neq j\;$ then $d(p(i,k),p(j,l))>4 D$ for
all $k\in\{1,\hdots,Q_{i}\}$, $l\in\{1,\hdots,Q_{j}\}$;\\
\item \label{iii} For all $i\in \{1,\hdots,s\}$ and for all
$k,l\in\{1,\hdots,Q_{i}\}$, there exists a sequence
$k_{1},k_{2},\hdots,k_{Q_{i}}$ of not necessarily distinct
elements of $\{1,\hdots,Q_{i}\}$ such that
$k=k_{1},\;l=k_{Q_{i}}$, and\\
$d\left(p(i,k_{j}),p(i,k_{j+1})\right)\leq 4 D$
for $j=1,\hdots,Q_{i}-1$.\\
\end{enumerate}

\noindent One notes that $\sum_{i=1}^{s} Q_{i}=Q$. We now define
for each $i=1,\hdots,s$,
$$f_{i}:  {\bf S}^{m}  \rightarrow  {\bf Q}_{Q_{i}}(Y)$$
such that $$\text{spt}(f_{i}(x))=\text{spt}(f(x))\cap
\left(\cup_{k=1}^{Q_{i}}{\bf B}(p(i,k),D)\right)$$ for each $x\in
{\bf S}^{m}$ where ${\bf B}(p(i,k),D)$ denotes the closed ball
with center $p(i,k)$ and radius $D$. We will now check that the
$f_i$ are well defined. For $i\in\{1,\hdots,s\}$ and $x\in {\bf
S}^{m}$ , there is at least $Q_{i}$ points in the support of
$f_{i}(x)$ since
$$\mathcal{S}\left(f(x),\sum_{l=1}^{s}\sum_{k=1}^{Q_{l}}\llbracket p(l,k)\rrbracket\right)=
\mathcal{S}(f(x),f(1,0,\hdots,0))\leq
\text{Lip}(f)\;|x-(1,0,\hdots,0)|\leq D.$$ Suppose that there
exists a point $p\in$
spt$\left(f_{i}(x)\right)\;\cap\;$spt$\left(f_{j}(x)\right)$ for
$i,j\in \{1,\hdots,s\}$ such that $i\neq j$. Then, there exist
$k\in\{1,\hdots,Q_{i}\},\; l\in\{1,\hdots, Q_{j}\}$ with
$d(p(i,k),p)\leq D$ and $d(p(j,l),p)\leq D$ hence
$d(p(i,k),p(j,l))\leq 2D$ which contradicts (\ref{ii}). By these
two observations, the $f_{i}$ are well defined and
$f=\sum_{i=1}^{s}f_{i}$. By construction, it is clear that
$\mathcal{S}(f(x),f(y))=\sum_{i=1}^{s}\mathcal{S}(f_{i}(x),f_{i}(y))$
for all $x,y\in {\bf S}^{m}$ hence Lip$(f_{i})\leq\;$Lip$(f)$ for
$i=1,\hdots,s$. By the definition of $f_{i}$ and (\ref{iii}), we
also notice that
$$\max\left\{\;d(p(i,1),p)\;:\;
p\in\text{ spt}(f_{i}(x))\right\}\leq 4D(Q_{i}-1)+D\leq D(4Q-3)$$
for $i=1,\hdots,s$ and $x\in {\bf S}^{m}$. Let $\{c_{xy}\}$ be a
$\gamma$-weakly convex geodesic bicombing on $Y$. We set
$$\begin{array}{l}
\theta :\mathbb{R}^{m+1}\backslash \{0\}\rightarrow {\bf S}^{m}
:x\mapsto \frac{x}{|x|}
\end{array}$$
and observe that Lip$(\theta|\{x\in\mathbb{R}^{m+1}\; : \;
|x|=r\})=1/r$ for all $0<r<\infty$. We can now define the
extension
$$\begin{array}{ll}
F:{\bf B}^{m+1}\rightarrow {\bf Q}_{Q}(Y)\\
\\
F(0)=\sum_{i=1}^{s}Q_{i}\llbracket p(i,1)\rrbracket,\\
\\
F(x)=\sum_{i=1}^{s}\sum_{j=1}^{Q_i}\left\llbracket
c_{p(i,1),q_j^i(x)}\left(|x|\;d\left(p(i,1),q_j^i(x)\right)\right)\right\rrbracket\\
\end{array}$$
where we denote $f_i\circ\theta(x)=\sum_{j=1}^{Q_{i}}\llbracket
q_{j}^{i}(x)\rrbracket$ for $i=1,\hdots,s$ and for each $0\neq
x\in {\bf B}^{m+1}$. We easily check that $F|_{{\bf S}^{m}}=f$.
Let $x,y\in {\bf B}^{m+1}$ such that $0<|x|\leq |y|$ and fix
$z=\frac{|x|y}{|y|}$.
Since $|x|=|z|$, we see that $|y|-|z|=|y|-|x|\leq |y-x|$ and it is
clear that $|x-z|\leq |x-y|$. It is also easy to check that
$\theta(z)=\theta(y)$. On the one hand, we compute\\
\begin{eqnarray*}
 & & \mathcal{S}(F(x),F(y))\\
 & \leq & \mathcal{S}(F(x),F(z))+\mathcal{S}(F(z),F(y))\\
  & = & \mathcal{S} \left(\sum_{i=1}^{s} \sum_{j=1}^{Q_i}\left\llbracket
c_{p(i,1),q_j^i(x)}\left(|x|\;d(p(i,1),q_j^i(x))\right)\right\rrbracket,
 \sum_{i=1}^{s}\sum_{j=1}^{Q_i}\left\llbracket
c_{p(i,1),q_j^i(z)}\left(|z|\;d(p(i,1),q_j^i(z))\right)\right\rrbracket \right)\\
  & & + \mathcal{S}\left(\sum_{i=1}^{s} \sum_{j=1}^{Q_i}\left\llbracket
c_{p(i,1),q_j^i(z)}\left(|z|\;d(p(i,1),q_j^i(z))\right)\right\rrbracket,\sum_{i=1}^{s}
\sum_{j=1}^{Q_i}\left\llbracket
c_{p(i,1),q_j^i(y)}\left(|y|\;d(p(i,1),q_j^i(y))\right)\right\rrbracket\right)\\
\end{eqnarray*}
where we can suppose that
$$\mathcal{S}\left(f_{i}\circ\theta(x),f_i\circ\theta(z)\right)=\mathcal{S}\left(\sum_{j=1}^{Q_{i}}\llbracket
q_j^i(x)\rrbracket,\sum_{j=1}^{Q_{i}}\left\llbracket
q_j^i(z)\right\rrbracket
\right)=\max_{j=1,\hdots,Q_{i}}d\left(q_j^i(x),q_j^i(z)\right)$$
for $i=1,\hdots,s$ and $q_j^i(z)=q_j^i(y)$ for $j=1,\hdots,Q_{i}$
and $i=1,\hdots,s$ since $\theta(z)=\theta(y)$. We conclude that
\begin{eqnarray*}
& & \mathcal{S}(F(x),F(y))\\
 & \leq & \max_{i=1,\hdots,s}\max_{j=1,\hdots,Q_{i}}
 d\left(c_{p(i,1),q_j^i(x)}\left(|x|\;d(p(i,1),q_j^i(x))\right),
 c_{p(i,1),q_j^i(z)}\left(|z|\;d(p(i,1),q_j^i(z))\right)\right)\\
 &  &  + \max_{i=1,\hdots,s}\max_{j=1,\hdots,Q_{i}}
 d\left(c_{p(i,1),q_j^i(z)}\left(|z|\;d(p(i,1),q_j^i(z))\right),c_{p(i,1),q_j^i(y)}\left(|y|\;d(p(i,1),q_j^i(y))\right)\right)\\
    & \leq & \gamma\;|x|\; \max_{i=1,\hdots,s}\max_{j=1,\hdots,Q_{i}} d\left(q_j^i(x),q_j^i(z)\right)
 + (|y|-|z|)\;\max_{i=1,\hdots,s}\max_{j=1,\hdots,Q_{i}}
 \;d\left(p(i,1),q_j^i(y)\right)\\
      & \leq & \gamma\; |x| \max_{i=1,\hdots,s} \mathcal{S}(f_{i}\circ\theta(x),f_i\circ\theta(z))
      +(|y|-|z|)\max_{i=1,\hdots,s}\max\left\{\;d(p(i,1),p)\;:\;
p\in\;\text{spt}(f_{i}\circ\theta(y))\right\}\\
        & \leq & \gamma\; |x-z|\max_{i=1,\hdots,s}\text{Lip}(f_{i})+D(4Q-3)|x-y|.\\
         & \leq & (\gamma+8Q-6)\text{Lip}(f)|x-y|.
\end{eqnarray*}
On the other hand, we compute
\begin{eqnarray*}
\mathcal{S}(F(x),F(0)) & = &
\mathcal{S}\left(\sum_{i=1}^{s}\sum_{j=1}^{Q_i}\left\llbracket
c_{p(i,1),q_j^i(x)}(|x|\;d(p(i,1),q_j^i(x))\right\rrbracket,\sum_{i=1}^{s}Q_{i}
\left\llbracket
p(i,1))\right\rrbracket \right)\\
   & \leq  & |x|\;\max_{i=1,\hdots,s}\max\{\;d(p(i,1),p)\;:\;
p\in\;\text{spt}(f_{i}\circ\theta(x))\}\\
       & \leq & (8Q-6)\text{Lip}(f)|x|.
\end{eqnarray*}
\qed

\section{Nagata dimension of ${\bf Q}_{Q}(Y)$}\label{nagata}

We begin by giving the precise definition of the Nagata dimension.
Suppose that $(Y,d)$ is a metric space and
$\mathcal{B}=(B_{i})_{i\in I}$ is a family of subsets of $Y$. The
family is called $D$-bounded, for some constant $D\geq 0$, if
diam$(B_{i})=\sup\left\{\;d(x,x')\;:\;x,x'\in B_{i}\right\}\leq D$
for all $i\in I$.
For $s>0$, the $s$-multiplicity of $\mathcal{B}$ is the infimum of
all $n\geq 0$ such that every subsets of $Y$ with diameter
$\leq\;s$ meets at
most  $n$ members of the family.\\

\begin{definition}
Let $Y$ be a metric space. The Nagata dimension dim$_{N}(Y)$ of
$Y$ is the infimum of all integers $n$ with the following
property: there exists a constant $c>0$ such that for all $s>0$,
$Y$ has a $cs$-bounded covering with $s$-multiplicity at most
$n+1$.
\end{definition}

\vspace{0.4cm}  The Section 2 in \cite{U. Lang T. Schlichenmaier}
gathers a number of basic properties of the Nagata dimension. We can
quote for instance the following. The topological dimension of a
metric space $Y$ never exceeds dim$_{N}(Y)$. The Nagata dimension is
a bilipschitz invariant and, as it turns out, even a quasisymmetry
invariant. The class of metric spaces with finite Nagata dimension
includes all doubling spaces, metric $(\mathbb{R}-)$trees, Euclidean
buildings, and homogeneous Hadamard manifolds, among others. Let us
study the Nagata dimension of ${\bf Q}_{Q}(Y)$ in accordance with
the Nagata
dimension of $Y$.\\

\begin{lem}\label{soit}
If $Y$ is a metric space with dim$_{N}(Y)=n<\infty$, then
dim$_{N}({\bf Q}_{Q}(Y))\leq (n+1)^{Q}-1$.
\end{lem}
\vspace{0.4cm} \noindent {\bf Proof.} Fix $s>0$. Since
dim$_{N}(Y)=n$, there exists $c>0$ such that $Y$ admits a
$cs$-bounded covering family $\mathcal{B}=(B_{i})_{i\in I}$ with
$s$-multiplicity at most $n+1$. For all $\sum_{j=1}^{Q}\llbracket
i_{j}\rrbracket\in {\bf Q}_{Q}(I)$, we fix
$B_{\sum_{j=1}^{Q}\llbracket i_{j}\rrbracket}:={\bf Q}_{Q}(Y)\cap
\left\{\sum_{j=1}^{Q}\llbracket x_{j}\rrbracket\; :\; x_{j}\in
B_{i_{j}}\text{ for }j=1,\hdots,Q\right\}$ and we define the
collection
$$\mathcal{B}^{*}=\left(B_{\sum_{j=1}^{Q}\llbracket i_{j}\rrbracket}\right)_{{\sum_{j=1}^{Q}\llbracket i_{j}\rrbracket}\in
{\bf Q}_{Q}(I)}.$$ One readily checks that $B^{*}$ is a
$cs$-bounded covering of ${\bf Q}_{Q}(Y)$. It remains to study the
$s$-multiplicity of $\mathcal{B}^{*}$. Let $A$ be a subset of
${\bf Q}_{Q}(Y)$ with diam$(A)\leq s$. It is clear that there
exist $A_{1},A_{2},\hdots, A_{Q}\subset Y$ such that
$A=\left\{\sum_{i=1}^{Q}\llbracket x_{i}\rrbracket\; :\;x_{i}\in
A_{i}\text{ for }i=1,\hdots,Q \right\}$ with diam$(A_{i})\leq s$
for $i=1,\hdots,Q$. We know that $A_{1},\hdots,A_{Q}$ meet
respectively at most $n+1$ members of $\mathcal{B}$. Consequently,
$A$
meets at most $(n+1)^{Q}$ members of $\mathcal{B}^{*}$.\qed\\

\begin{coro}\label{aline}
If $Y$ is a Banach space with a finite algebraic dimension, then
the Nagata dimension of ${\bf Q}_{Q}(Y)$ is finite.
\end{coro}
\vspace{0.4cm} \noindent {\bf Proof.} If $Y$ has a finite
algebraic dimension, the unit ball of $Y$ is precompact.
Consequently, $Y$ is doubling and has a finite Nagata dimension
hence dim$_{N}({\bf
Q}_Q(Y))<\infty$ by Lemma \ref{soit}.\qed\\

We are now able to prove the results mentioned in the
introduction. On the one hand, Theorem \ref{exto} immediately
 ensues from Lemma \ref{complete}, Proposition
\ref{connexe} and Theorem 1.5 in \cite{U. Lang T. Schlichenmaier}.
On the other hand, Theorem \ref{final} is an immediate consequence
of Lemma \ref{complete}, Proposition \ref{connexe}, Lemma
\ref{soit}, Corollary \ref{aline} and Corollary 1.8 in \cite{U.
Lang T. Schlichenmaier}.

\section{Acknowledgements}
The author is thankful to Thierry De Pauw for suggesting the
problem and helpful discussions.

\end{document}